\documentclass{gtart}
\usepackage{amsmath,amssymb}
 
\input gtmonout
\volumenumber{2}
\volumeyear{1999}
\volumename{Proceedings of the Kirbyfest}
\pagenumbers{321}{333}
\papernumber{17}
\received{20 January 1999}\revised{10 October 1999}
\published{20 November 1999}

\newtheorem{thm}{Theorem}[section]
\newtheorem{lm}[thm]{Lemma}

\theoremstyle{definition}
\newtheorem{example}[thm]{Example}
\newtheorem{df}[thm]{Definition}

\input epsf

\def \r {\text{\em Res}}

\def \Q {{\mathbf{Q}}}

\begin{document}
 
\title{Fermat limit and congruence of 
Ohtsuki invariants}

\author{Xiao-Song Lin\\Zhenghan Wang}
\address{Department of Mathematics, University of California\\ 
Riverside, CA 92521, USA\\\smallskip
\\Department of Mathematics, Indiana University\\Bloomington, 
IN 47405, USA}
\email{xl@math.ucr.edu, zhewang@indiana.edu}
\asciiaddress{Department of Mathematics, University of California\\ 
Riverside, CA 92521, USA\\
\\Department of Mathematics, Indiana University\\Bloomington, 
IN 47405, USA}

\begin{abstract}
By calculating the Fermat limit of certain $q$--Fermat functions,
we get explicit surgery formulae for the second and third Ohtsuki
invariants 
for homology 3--spheres. The surgery formula of the second Ohtsuki
invariant $\lambda_2$, 
combined 
with an argument using the 
general theory of finite type invariants of homology 3--spheres, leads to a
congruence relation between the first and second Ohtsuki invariants:
$\lambda_1= 2\lambda_2$ mod 24.
\end{abstract}

\asciiabstract{By calculating the Fermat limit of certain q-Fermat
functions, we get explicit surgery formulae for the second and third
Ohtsuki invariants for homology 3-spheres. The surgery formula of the
second Ohtsuki invariant \lambda_2, combined with an argument using
the general theory of finite type invariants of homology 3-spheres,
leads to a congruence relation between the first and second Ohtsuki
invariants: \lambda_1=2\lambda_2 mod 24.}



\maketitle

\cl{\small\it Dedicated to Rob Kirby on the occasion of his 60th birthday} 

\section{Introduction} In \cite{O1}, Ohtsuki extracted a series of
rational 
topological invariants denoted by 
$\lambda_1,\lambda_2,\dots,\lambda_n,\dots$ of 
oriented homology 3--spheres from the $SO(3)$ quantum invariants of 
Reshetikhin and Turaev \cite{RT,KM}. 
Physically, they correspond to the coefficients of the 
asymptotic expansion of Witten's Chern--Simons path integral at 
the trivial connection as shown by Rozansky \cite{R}. This construction 
stimulated a series of works on finite type invariants of homology 
3--spheres.

The purpose of this article is to present the following theorem about
a congruence relation between the first and second Ohtsuki invariants:

\begin{thm} For every oriented homology 3--sphere $M$, we have
\begin{equation}\label{cong}
\lambda_1(M)=2\lambda_2(M) \qquad
\text{\rm mod 24}.
\end{equation}
\end{thm}

\noindent{\bf Remark}\qua It is known that $\lambda_1(M)/6$ is equal to the 
Casson invariant of $M$ \cite{Mu}. Here the Casson invariant $\lambda(M)$
is 
normalized such that
$\lambda=1$ for the Poincar\'e homology 3--sphere $\Sigma(2,3,5)$
($+1$--surgery 
on the right trefoil knot). It is also known \cite{LiW} that
$\lambda_2(M)\in
3{\mathbb Z}$ for every $M$ (a simpler proof will be given here in Section
5).
So the congruence (\ref{cong}) could also be 
expressed as
\begin{equation}\label{Cong}
\frac{\lambda_1(M)}{6}=\frac{\lambda_2(M)}{3} \qquad
\text{\rm mod 4}.
\end{equation}

The proof of the congruence (\ref{cong}) goes roughly as follows: 
In \cite{LiW},
we introduced the notion of Fermat functions and Fermat limit. By an 
explicit calculation of the Fermat limit of a certain $q$--Fermat function 
involving weighted Gauss sums, we get a surgery formula for $\lambda_2$.
The 
congruence (1) is then first proved using the surgery formula in the case
when
$M$ is obtained from Dehn surgery on a knot in $S^3$. 

Explicit surgery formulae for $\lambda_n$ were also given in \cite{LiW}.
These
formulae were used to prove that $\lambda_n$ is a finite type invariant of 
order $3n$ in \cite{KrSp}. See also \cite{O3} and \cite{L,L2}. 
We then use the general 
machinery of finite type invariants of homology 3--spheres to argue that
the
congruence (1) for homology 3--sphere obtained from Dehn surgery on knots 
implies that it holds for all homology 3--spheres.

We originally expressed our congruence in the form (2). When Rob Kirby saw
it, 
he suggested that we put it in the form (1) because of the significance of
24 
in the study of 3--manifolds. He also asked whether there are general 
congruence 
relations among $\lambda_1,\lambda_2,\dots,\lambda_n$ similar to (1). We
are 
unable to have any hint about his question even for $\lambda_3$ because
the 
formulae are rather involved.  

It is possible that the congruence (1) is related with the connection 
between modular forms and quantum invariants of 3--manifolds suggested in
\cite{LZ}. 

The article is organized as follows: We will outline some calculation in
\cite{LiW} leading to a surgery formula of $\lambda_n$ in Sections 2, 3,
4. 
The reader will notice that Formula (3) is better expressed here than 
in \cite{LiW}, and Formula (6) (the surgery formula for $\lambda_3$) 
is new here. In Section 5 we will 
prove Theorem 1.1 using surgery formulae for $\lambda_{1,2}$ and the
theory 
of finite type invariant of homology 
3--spheres.

\rk{Acknowledgements}The first author is supported in part by an NSF grant 
and the second author is supported by an NSF postdoctoral fellowship.

\section{Fermat Functions and Their Residues} 

We fix some notation that will be in force throughout this 
article.  Let ${\mathbb P}$ denote the set of all {\it odd} prime numbers,
and ${\mathbb Q}$ the rational numbers.

Given $r\in {\mathbb P}$, 
\begin{itemize}
\item ${\mathbb Z}_r$ is the ring of $r$--adic integers.
\item We use $\bar{m}$ to denote the integer in $\{1,\cdots, r-1\}$
such that $m\cdot \bar{m}=1$ mod $r$.
\item We will use $q$ to denote the $r$--th root of unity, 
$e^{2 \pi\sqrt{-1}/r}$.  
The quantum integer $[k]_q$ is
$$[k]_{q}=\frac{q^{\frac{k}2}-q^{-\frac{k}2}}{q^{\frac12}-q^{-\frac12}}.$$ 
\item For $R=\mathbb Z$ or ${\mathbb Z}_r$, $O((q-1)^k;R)$ 
stands for a complex number 
of the form $u(q-1)^k$ for some $u\in R[q]$.
\item The {\it Gauss sum} is $G_0(q)=\sum_{k=0}^{r-1}q^{k^2}$. The {\it
weighted Gauss sum} is\break $G_{2l}(q)=\sum_{k=0}^{r-1}k^{2l}q^{k^2}$.
\end{itemize}

\begin{df}
Suppose  $f$ is a function $f\co  {\mathbb P} \rightarrow {\mathbb Q}$ 
such that $f(r)\in{\mathbb Z}_r$ for all sufficiently large $r$, 
$f(r)$ is a {\em Fermat function} if there exists a 
rational number $\lambda=m/n$ (in simplified form) 
independent of $r$ such that
when $f(r)=m'/n'$ (in simplified form), then $m'n= mn'$ mod $r$ for all 
sufficiently large $r$.    

If $f$ is a Fermat function, we will call the rational number 
$\lambda$ the {\it residue} of $f$ and denote it by $\text{Res}(f)$. 
\end{df}

\begin{lm} The residue of a Fermat function is unique.
\end{lm}

\begin{proof} If there is another rational number $\lambda'$ such that 
$f(r)=\lambda'$ mod $r$ for all 
 sufficiently large primes $r$, then the numerator of 
the rational number $\lambda-\lambda'$ will be divisible 
by all sufficiently large primes $r$. 
Therefore, $\lambda-\lambda'=0$,  ie $\lambda=\lambda'$.
\end{proof}

A typical Fermat function comes from Fermat's little theorem (and hence
the 
name). For any rational $a\neq 0$,  
the function $f(r)=a^{r-1}$ is a Fermat function with $\text{Res}(f)=1$. 
There are many other Fermat functions. For example, the function
$f(r)=(r-1)/2$
is a Fermat function whose residue is $-1/2$. Also, the function 
$f(r)=(r-1)!$ is a 
Fermat function because of Wilson's theorem, which says that 
$(r-1)!=-1$ mod $r$ \cite{HW}. On the other hand, an example of a
non-Fermat function is given by $f(r)=(\frac{r-1}{2})!$.  See \cite{HW}
for a
discussion of the residue of $(\frac{r-1}{2})!$, which turns out to depend
on 
$r$ in a quite complicated way. 

\begin{lm}
The set of all Fermat functions is a ring over the rationals ${\mathbb
Q}$.
\end{lm}

This follows easily from the following facts.
\begin{lm} Suppose that $f$ and $g$ are Fermat functions. Then
\begin{enumerate}
\item for rational numbers $\alpha$ and $\beta$, $\alpha f+\beta g$ is a 
Fermat function whose residue is $\alpha \cdot \r (f)+\beta \cdot \r (g)$;
\item $f\cdot g$ is a Fermat function whose residue is $\r (f)\cdot\r
(g)$;
\item if $\r (g)\neq 0$, then $f/g$ is a Fermat function with residue
$\r (f)/\r (g)$.
\end{enumerate}
\end{lm} 

Thus, in particular, every polynomial function of $r$ with rational 
coefficients is a Fermat function whose residue is its constant term. 
And every rational function of $r$ with rational coefficients 
is a Fermat function if the constant term of the denominator 
is not zero. Its residue is the value of the function at $r=0$.
More generally, we have 

\begin{lm}
Suppose $F(r)=\frac{P(r)}{Q(r)}$ is a rational polynomial of $r$, and
$f(r)$ 
is a Fermat function with residue $\lambda$.  If $Q(\lambda)\neq 0$, 
then the composite function $F(f(r))$ is a Fermat function.
\end{lm} 

We give another example of Fermat functions and the 
details can be found in \cite{O1,LiW}.

\begin{example} For a fixed integer $k$, the function
$$D_k(r)=\frac{\left(\frac{r-1}2\right)!}{\left(\frac{r-1}2-k\right)!}
$$
is a Fermat function whose residue is given by
$$\text{Res}(D_k)=\begin{cases} 
(-\frac12)(-\frac12-1) \cdots (-\frac12-(k-1))
& \text{for $k>0$,}\\
1  & \text{for $k=0$,}\\
\frac{1}{(-\frac12+1)\cdots (-\frac12-k)} &\text{for $k<0$.}
\end{cases}$$
\end{example}

Recall that $f(r)=(\frac{r-1}{2})!$ is not a Fermat function.

\section{$q$--Fermat Functions and Fermat Limit}

The $SO(3)$ quantum invariants of 3--manifolds take values in ${\mathbb
Z}[q]$ 
when appropriately normalized.  We like to study the so-called
perturbative 
expansion of these invariants.  This leads us to think of  
$q$ as a variable, and study the expansion around $q=1$.  
In this section, we describe the $q$--analogue of Fermat functions when 
the function takes values in ${\mathbb Z}_r$ adjoint with the $r$-th 
root of unity $q$, ${\mathbb Z}_r[q]$.   
It is reasonable to expect that the corresponding object of 
residue lies in the ring ${\mathbb Q}[[t-1]]$ 
of formal power series in $t-1$ with coefficients in $\mathbb Q$, where
$t$
is a formal variable in place of $q$.
  
Let $c(r)$ be an integer sequence indexed by $r$ with
$$\lim_{r\rightarrow +\infty} c(r)=+\infty\qquad\text{and}
\qquad c(r)\leq r-2.$$
Given a complex function $f_{q}(r)$  on $\mathbb P$ which takes values 
in ${\mathbb Z}_r[q]$.  
Fix a non-negative integer $n$. For sufficiently large $r$, we write
$$
\begin{aligned}
f_q(r)=& a_{r,0}+a_{r,1}(q-1)+\cdots+a_{r,n}(q-1)^n+\\
&\cdots+a_{r,c(r)}(q-1)^{c(r)}+O((q-1)^{c(r)+1},{\mathbb Z}_r)
\end{aligned}
$$ 
for some $a_{r,n}\in {\mathbb Z}_r$.
 
Notice that although
the expansion $f_q(r)$ is not unique as $q$ is not a free variable,
$a_{r,n}\in{\mathbb Z}/r{\mathbb Z}$ is well defined for $0\leq n\leq
c(r)$.
This is because of $r=O((q-1)^{r-1};\mathbb Z)$, and the relation
in ${\mathbb Z}_r[q]$: $\sum_{i=0}^{r-1} q^{i}=0$. Thus, if $a_{r,n}$
is a Fermat function, its residue depends only on the 
function $f_q(r)$.

\begin{df} The complex function $f_q(r)$  has
a {\em Fermat limit}, denoted by $\text{f-lim}\,f_q(r)$,
if each $a_{r,n}$, thought as a function of $r$,
is a Fermat function with $\r(a_{r,n})=\lambda_n$. By definition,
$$\text{f-lim}\,f_q(r)=\sum_{n=0}^{\infty}\lambda_n(t-1)^n
\in {\mathbb Q}[[t-1]].$$
We will call $f_q(r)$ a {\em $q$--Fermat function} if its Fermat limit
exists.
\end{df}

The Fermat limit of $f_q(r)$ is well-defined and 
if it exists, it is unique.  
The uniqueness of a Fermat limit follows from the uniqueness of residues 
of Fermat functions (Lemma 2.1).

The following are some basic properties of the Fermat limit.

\begin{itemize}
\item The $q$--analogues of both Lemma 2.2 and Lemma 2.3 hold.
\item If $f(r)$ is a Fermat function which takes values in $\mathbb Z$
with 
residue $\lambda$, 
then $q^{f(r)}$ is a $q$--Fermat function whose Fermat limit is 
$t^{\lambda}$, where it is understood that $t^{\lambda}$ is 
expanded as a power series of $t-1$.
\end{itemize}

Let $$\tilde G_{2l}(q)=(q-1)^l\frac{G_{2l}(q)}
{G_{0}(q)}.$$
The following theorem (see \cite{LiW} Lemmas 6.1 and 6.2) 
is very important for our study of the 
$SO(3)$ quantum invariants of 3--manifolds.

\begin{thm} $\tilde G_{2l}(q)$ is a $q$--Fermat function and 
\begin{equation}
\text{\em f-lim}\,\tilde G_{2l}(q)
=A_{l}\cdot \left(\frac{t-1}{\log t}\right)^{l}
\end{equation}
where $A_{l}=(-\frac{1}{2})^{l}(2l-1)!!=(-\frac{1}{2})^{l}\cdot 1\cdot 3
\cdots (2l-1),\; and \;\;A_0=1$.
\end{thm}

One may compare (3) with
$$\frac{ \int_{-\infty}^{\infty} x^{2l}e^{-x^{2}}dx}
{\int_{-\infty}^{\infty}e^{-x^{2}}dx} =(-1)^{l}A_{l}.$$

\section{Surgery Formulae of Ohtsuki invariants}

Let $L=K_{1}\cup K_{2}\cup\cdots \cup K_{\# L}$ be an oriented link in
$S^{3}$, where we use $\# L$ to denote the number of components of $L$.  
We will denote the unknot by $O$, and the empty link by
$\emptyset$.
  
The Conway polynomial $\nabla (L;z)\in{\mathbb Z}[z]$ is defined
by 
$$\begin{cases} \nabla (O;z)=1 \\
\nabla(\emptyset;z)=0 \\
\nabla (L_{+};z)-\nabla (L_{-};z)=-z\cdot \nabla(L_{0};z).
\end{cases}
$$ 
Here, as usual, $L_{+}$, $L_{-}$ and $L_{0}$ are the links which have
plane
projections identical to each other except 
in one small disk where their projections are a positive crossing, a
negative crossing and an orientation preserving smoothing of that
crossing, respectively. 
Note the negative sign on the right hand side of the skein relation. 
Its effect is to change the usual Conway polynomial by a normalization
factor
$(-1)^{\# L-1}$. So for a knot $K$, $\nabla(K;z)$ is the same as the usual 
Conway polynomial of $K$. In particular 
$$\nabla(K;z)=1+c_2(K)z^2+c_4(K)z^4+\cdots+c_{2k}(K)z^{2k}.$$
The Jones polynomial $V(L;t)\in
{\mathbb Z}[t^{\frac12}, t^{-\frac12}]$ is defined by 
$$\begin{cases}
V(O;t)=1;\\
V(\emptyset;t)=(t^{\frac12}+t^{-\frac12})^{-1}\\
tV(L_{+};t)-t^{-1}V(L_{-};t)=(t^{\frac12}-t^{-\frac12})V(L_{0};t).\end{cases}
$$
Note that our normalization differs
from the usual definition of the Jones polynomial.  Actually it
is obtained from the usual Jones polynomial by changing $t$ to
$t^{-1}$ and multiplying by $(-1)^{\# L-1}$.  We put
$$X(L;t)=\frac{V(L;t)}{(t^{\frac{1}{2}}+t^{\frac{-1}{2}})^{\# L -1}}.$$
Then
$$X(O;t)=X(\emptyset:t)=1.$$
We also put
$$\Phi(L;t)=\sum_{L'\subset
L}(-1)^{\# L-\# L'} X(L';t)$$
where the sum runs over all sublinks of $L$
including the empty link and $L$ itself.  
We have $\Phi(O;t)=0$ and define
$\Phi(\emptyset;t)=0$.  

Further, we put 
$$\Phi_{i}(L)=\left.\frac{d^{i}\Phi(L;t)}{dt^{i}}\right|_{t=1}$$
so that
$$\Phi(L;t)=\sum_{i=0}^\infty\frac{\Phi_i(L)}{i!}(t-1)^i.$$
Finally for each $i\geq 1$, we set
$$\phi_{i}(L)=\frac{(-2)^{\# L}}{(\# L+i)!}\cdot \Phi_{\# L +i}(L).$$
The invariants $\phi_i$ are the basic link invariants we use
to express $\lambda_{n}$.

A link is called an {\em algebraically split link} (ASL) if the 
linking number between every pair of components is 0. In particular,
a knot is always an ASL.
A framed ASL is said to be {\em unit framed} if the framing of each
component
is $\pm 1$.  Given a link $L$ and a positive integer $m$, we use
$L^{m}$ to denote the 0--framed $m$--parallel of $L$, ie, each component in
$L$ is replaced by $m$ parallel copies having linking number zero with
each
other.  So
$L^{m}$ is  an ASL if $L$ is.  
When $L$ is ordered, sublinks of $L^m$ will be in one--one
correspondence with $\mu$--tuples $(i_1,\dots,i_\mu)$, where $\mu=\#L$,
in such a way that $L'$ has $i_\xi$ parallel copies of the $\xi$-th
component of
$L$, $0\leq i_\xi\leq m$. 
If $L$ is a framed link, $L^m$ and all its sublinks will inherit a framing
from
$L$. If the $\xi$-th component of $L$ is framed by $f_\xi$,
$\xi=1,\dots,\mu$, we denote
$$f_L=\prod_{\xi=1}^{\mu}f_\xi.$$
The existence of Ohtsuki invariants is based on two key ingredients.
The first one is the following property of Jones polynomial for cablings:

\begin{lm}{\em (Ohtsuki, Proposition 3.4 in \cite{O1})} As a power series
in $t-1$, 
$$ \Phi(L^{{\bf i}};t)=O((t-1)^{|{\bf i}|+\text{\em max}({\bf
i})})\in\Q[[t-1]]$$
if $L$ is an ASL.  Here if ${\bf i}=(i_1, \cdots , i_{\sharp L}),$
then $|{\bf i}|=i_1+\cdots +i_{\sharp L}$, and 
{\em max({\bf i})} is the largest of $i_j$.
\end{lm}

This implies that the following formal power series is 
well-defined for every ASL:
$$\tilde{\Phi}(L)=\sum_{l=0}^{\infty}\sum_{L'\subset L^{l}\backslash
L^{l-1}}
\frac{\Phi(L;t)}{(t-1)^{\sharp L'}}
$$
Another key ingredient is the existence of the Fermat limit 
of the following function:
$$ 
\begin{aligned}
H_{i, f}(q)&=(-f)\cdot \left(\frac{f}{r}\right)
\cdot q^{3\cdot {\bar{4}}\cdot f-{\bar{2}}}\cdot 
\frac{(q-1)^{i+1}}{G_{0}(q)}\\  
&\cdot \left\{ \sum_{k=1}^{\frac{r-1}{2}}
q^{{\bar{4}}\cdot f(k^{2}-1)} [k]_q\sum_{j=0}^{\frac{k-1}{2}} (-1)^{j} 
\binom{k-j-1}{j}\binom{k-2j-1}{i}
[2]_q^{k-2j-1} \right\},
\end{aligned}$$
\noindent where $G_{0}(q)$ is the Gauss sum
and $\left(\frac{f}{r}\right)$ is the Legendre
symbol. Note that $H_{i, f}(q)\in {\mathbb Z}_r[q]$. 

The existence of the Fermat limit of $H_{i,f}(q)$ is established by 
Ohtsuki \cite{O1}, 
and our computation of this limit is the key 
to the explicit formulae of Ohtsuki invariants.
Our computation relies on Theorem 3.2.

Suppose $L$ is a unit-framed ASL, then the $SO(3)$ quantum invariant 
of the oriented homology 3--sphere $S^{3}_{L}$ can be expressed as
follows:
$$
\tau_{r}(S^{3}_{L})=
\sum_{l=0}^{\frac{r-3}{2}}\sum_{L'\subset L^{l}\backslash L^{l-1}}
\frac{\Phi(L';q)}{(q-1)^{\sharp L'}}\cdot H_{{\bf i(L')},{\bf f}}(q),$$
here ${\bf i(L')}$ is the corresponding $\sharp L$-tuple of $L'$.
For $\mu$-tuples ${\bf i}=(i_1, i_2,\cdots , i_{\mu})$, and 
 ${\bf f}=(f_1, f_2,\cdots , f_{\mu})$,  
$H_{{\bf i}, {\bf f}}=\prod_{j=1}^{\mu}H_{i_j,f_j}$.

Now let $\text{f-lim}\,H_{i,f}(q)=H_{i,f}(t)$. Then as a power series of 
$t-1$, we have 
$$ 
\begin{aligned}
\text{f-lim}\,\tau_{r}(S^{3}_{L}) 
&=\sum_{l=0}^{\infty}\sum_{L'\subset L^{l}\backslash L^{l-1}}
\frac{\Phi(L';t)}{(t-1)^{\sharp L'}}\cdot H_{{\bf i(L')},{\bf f}}(t)\\
&=1+\lambda_{1}(t-1)+\lambda_{2}(t-1)^{2}+\cdots  
\end{aligned} $$
Explicit surgery formulae for $\lambda_{i},i=1,2,3$ have been 
worked out using the above formulation. We need the following calculation:
\begin{itemize}
\item $H_{0,f}(t)=1$
\item $H_{1,f}(t)=(-f)\cdot (t+1)$
\item $H_{2,f}(t)=f\cdot \frac{(t+1)^{2}}{2}
\cdot [1+f+4\sum_{m=1}^{\infty}g_{1,m}(t-1)^{m-1}]$
\item $H_{3,f}(t)=(-f)\cdot \frac{(t+1)^{3}}{6}\cdot t^{\frac{1+f}{2}}
\cdot [1+24 \cdot \sum_{m=0}^{\infty}g_{1,m+3}(t-1)^{m+1}],$
\end{itemize}
where $g_{i,j}$ are constants determined by 
$$\sum_{m=0}^{\infty}g_{l,m}(t-1)^{m}=A_{l}\cdot\left(\frac{t-1}{\log t}
\right)^{l},\qquad\qquad l\geq 0$$ 
(both sides are power series of $t-1$).

\begin{thm} Let $L$ be a unit framed ASL and $S^{3}_{L}$ be the
homology 3--sphere obtained from Dehn surgery on $L$. Let $\lambda_i$ 
 be the $i$-th Ohtsuki invariant. Then
\begin{equation}
\lambda_{1} (S^{3}_{L})=\sum_{L'\subset L}f_{L'}\,\phi_{1}(L')
\end{equation}
and
\begin{equation}
\lambda_{2} (S^{3}_{L})=\sum_{L'\subset L}f_{L'}\, \phi_{1}(L') \,
 \frac{ \# L'}{2} 
+\sum_{L'\subset L^{2}}f_{L'}\, \phi_{2}(L')\, 
\frac{1}{2^{s_{2}(L')}}
\end{equation}
and 
\begin{eqnarray}
\lambda_{3} (S^{3}_{L}) =
\sum_{L'\subset L}f_{L'}\, \phi_{1}(L') \,
  \frac{ \# L'(\# L' -1)}{8}\hspace{1in} \mbox{}\nonumber \\  
 +\sum_{L'\subset L^{2}}f_{L'}\, \phi_{2}(L')\, 
 \frac{s_{1}+2s_{2} 
+\frac{1}{3} \sum_{i_{\xi}=2}f_{\xi}}{2^{s_{2}+1}}
+\sum_{L'\subset L^{3}}f_{L'}\, \phi_{3}(L')\, 
 \frac{1}{2^{s_{2}+s_{3}}\cdot 3^{s_{3}}}.
\end{eqnarray}
Here, if $L'$ corresponds to the $\mu$--tuple $(i_1,\dots,i_\mu)$,
$s_{j}(L')=
\,\# \{i_\xi\,;\,i_\xi=j\}$.
\end{thm}

See Theorem 5.1 in \cite{LiW}. The formula (6) is new here. It is proved 
following the same but more tedious calculation used in the proof of (4)
and (5).
So we will not repeat it here.  

In the rest of this section, we describe the invariants
$\lambda_1, \lambda_2$ for homology 3--spheres obtained from 
$1/n$--Dehn surgery on knots.

Let $c_{4}$ be the coefficient of $z^{4}$ in the Conway
polynomial of $K$ and $v_i$ be the $i$-th derivative of $V(K;e^h)$ at
$h=0$.

\begin{thm}
Let $n$ be an integer, and  
$S^{3}_{K,1/n}$ be the homology 3--sphere obtained from $1/n$--Dehn
surgery on a knot $K$.  Then
\begin{equation}
\lambda_{1}(S^{3}_{K,1/n})=-n\cdot v_{2}(K)
\end{equation}
and
\begin{equation}
\lambda_{2}(S^{3}_{K,1/n})=\frac{n}{2}\, v_{2}(K)
-\frac{n}{3}\, v_{3}(K)+n^{2}\, [v_{2}(K) +\frac{5}{3}v_{2}^{2}(K)
-60c_{4}(K)].
\end{equation}
\end{thm}

Some examples are as follows.

\begin{example}
(i) If $M$ is the Poincar\'e homology 3--sphere $\Sigma(2,3,5)$
($+1$--surgery 
on the right-handed trefoil knot), then 
$\lambda_1(M)=6$ and $\lambda_{2}(M)=39$.

(ii) If $M$ is the homology 3--sphere 
$\Sigma(2,3,7)$ ($+1$--surgery on the 
left-handed trefoil knot), then 
$\lambda_1(M)=6$ and $\lambda_{2}(M)=63$.
\end{example}

We may check that the congruence (1) or (2) holds in these examples.

\section{Proof of Theorem 1.1}

We first consider the special case when $M=S^3_{K,1/n}$. 

\begin{lm} For every knot $K$, $v_2(K)\in6{\mathbb Z}$ and $v_3(K)\in36
{\mathbb Z}$.
\end{lm}

\begin{proof} Using the fact that the knot invariants $v_2$ and $v_3$ are
of finite type (see eg \cite{BN}), we only need to check the desired
congruences on the (left-handed) trefoil knot (for $v_2$) as well as
on the figure-8 knot (for $v_3$). We have 
\begin{itemize}
\item $V(\text{terfoil};e^h)=e^h-e^{4h}+e^{3h};$
\item $V(\text{figure-8};e^h)=1+e^{-2h}+e^{2h}-e^{-h}-e^{h}.$
\end{itemize}

Thus, $v_2(\text{trefoil})=1-4^2+3^2=-6$,
$v_3(\text{trefoil})=1-4^3+4^3=-36$,
and $v_3(\text{figure-8})=(-2)^3+2^3-1^3-(-1)^3=0$.
\end{proof}

By Lemma 5.1 and Equations (7, 8), we see $\lambda_1(M)=2\lambda_2(M)$ mod 
24 if $M=S^3_{K,1/n}$. We also see $\lambda_2\in3{\mathbb Z}$ here. 

If $M=M_1\#M_2$, we have
$$\begin{aligned}
\,&\lambda_1(M)=\lambda_1(M_1)+\lambda_1(M_2),\\
&\lambda_2(M)=\lambda_2(M_1)+\lambda_2(M_2)+\lambda_1(M_1)\lambda_1(M_2).
\end{aligned}
$$
So, the desired congruences hold for connected sums of homology 3--spheres
obtained from surgery on knots.

We consider the general case now.

There is a parallel theory of finite type invariants for 
oriented homology 3--spheres (see \cite{O2,O3,Li}) and Ohtsuki invariants 
$\lambda_n$ are known to be finite type invariants of order $3n$. 
In particular, for any ASL with $\sharp L> 6$, we have 
$$\sum_{L'\subset L}(-1)^{\sharp L'} \lambda_{1,2}(S^{3}_{L'})=0.$$
So essentially, we only need to calculate $\lambda_{1,2}$ and check the 
congruences on the homology 3--sphere $M_*$
obtained from $+1$ surgery on a 6--component link $L$ coming
from the graph {\epsfxsize=0.15 truein \epsfbox{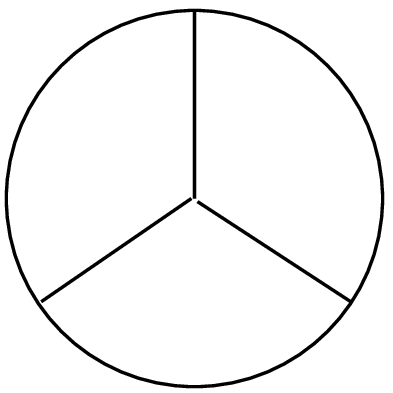}} in the 
following way: First replace each edge with a circle, and then hook the 
circles together like Borromean rings for each tri-valent vertex.

A direct calculation of $\lambda_2$ on $M_*$ involves the 
computation of the Jones polynomial of a 12--component link with 96
crossings,
a task too tedious to pursue.

Instead, let us first blow down four circles in $L$ to get a two component
link $J\cup K$. Both components $J$ and $K$ are now knotted so that we
cannot 
blow down further. But the link $J\cup K$ has the following nice property:

There are collections of circles $\{C_1,C_2,C_3,C_4\}$ (each $C_i$
contains two circles) in the complement of 
$J\cup K$ so that
\begin{enumerate}
\item $C_1\cup C_2\cup C_3\cup C_4$ is a trivial link;
\item the link $J\cup K\cup C_1\cup C_2\cup C_3\cup C_4$ is an ASL and 
is surgery 
equivalent to a separate link $(J\cup C_1\cup C_2)\amalg(K\cup C_3\cup
C_4)$;  
\item blowing down all circles in any non-empty subset of
$\{C_1,C_2,C_3,C_4\}$
will unknot either $J$ or $K$. 
\end{enumerate}

Therefore, the values of $\lambda_{1,2}$ on $M_*$ are equal to
a linear combination (with integer coefficients) of that of
$\lambda_{1,2}$
on several homology 3--spheres obtained from surgery on knots, 
and a connected sum of such homology 3--spheres, respectively. 
Thus, $\lambda_2\in3{\mathbb Z}$ and $\lambda_1=2\lambda_2$ mod 24 also
hold. 
This finishes the proof of Theorem 1.1 as well as the fact that
$\lambda_2\in
3{\mathbb Z}$ for all homology 3--spheres.

The existence of such collections of circles $\{C_1,C_2,C_3,C_4\}$ will
be quite obvious once we have a good picture of $J\cup K$. It is a good 
exercise in the so-called Kirby calculus, which is a wonderful tool Rob
Kirby 
invented for our benefit. And we will leave the exercise to those who 
has not yet had enough chance to play with this tool.

\Addresses\recd

\end{document}